\documentclass[11pt]{article}
\usepackage{amsmath}
\usepackage[all]{xy}
\usepackage{amsfonts}
\usepackage{amssymb}
\usepackage{amscd}
\usepackage{amsthm}
\usepackage{latexsym}
\usepackage{amsbsy}
\newtheorem{lem}{Lemma}[section]
\newtheorem{prop}{Proposition}[section]
\newtheorem{theorem}{Theorem}[section]
\newtheorem{corollary}{Corollary}[section]
\newtheorem{definition}{Definition}[section]
\newtheorem{example}{Example}[section]
\newtheorem{remark}{Remark}

\newcommand{\ten}{\otimes_R}

\newcommand{\ot}{\otimes}

\newcommand{\ts}{\widetilde{S}}

\newcommand{\sw}[1]{^{(#1)}}
\newcommand{\Hc}{\mathcal{H}}
\newcommand{\cp}{\mathcal{P}}
\date{Received 27 July 2004}
\begin{document}
\title
{Para-Hopf algebroids and their cyclic cohomology}

\author {M. Khalkhali,~~~ B. Rangipour,
\\\texttt{~masoud@uwo.ca ~~~~brangipo@uwo.ca}
\\Department of Mathematics
\\University of Western Ontario }
\maketitle
\begin{ }
\noindent {\bf Abstract.}
We introduce the concept of {\it para-Hopf algebroid} and
define  their cyclic cohomology in the spirit of Connes-Moscovici
cyclic cohomology for Hopf algebras. Para-Hopf algebroids are
closely related to, but different from, Hopf algebroids. Their
definition is motivated by attempting to define a  cyclic
cohomology theory for Hopf algebroids in general. We show that
many of Hopf algebraic  structures, including the
Connes-Moscovici algebra $\mathcal{H}_{FM}$, are para-Hopf algebroids.\\



\end{ }


\section{Introduction}
Broadly speaking, Hopf algebroids are quantizations of groupoids. More precisely, they are
 the noncommutative and  non-cocommutative algebraic analogues of Lie groupoids and Lie algebroids.
In their study of  the index theory of transversely
elliptic operators and in order to obtain more canonical formulas which
work for non-flat transversals, Connes and Moscovici \cite{achm01} had to replace
their Hopf algebra $\mathcal{H}_n$
 by a so called extended Hopf algebra of {\it transverse differential operators}
$\mathcal{H}_{FM}$. Here $FM$ denotes the frame bundle of an
$n$-dimensional manifold $M$.  In fact $\mathcal{H}_{FM}$ is a
bialgebroid in the sense of \cite{ lu, Xu} and has a so called twisted
antipode $\tilde{S}$ with $\tilde{S}^2=1$; precise definitions will
be given in Section 2. To define a cyclic cohomology theory   for $\mathcal{H}_{FM}$,  Connes and
Moscovici used the natural  action of
$\mathcal{H}_{FM}$ on the  algebra $\mathcal{A}_{FM}$ associated
to the frame bundle of $M$ and defined the cocyclic module of
$\mathcal{H}_{FM}$ as a certain submodule of the cocyclic module
of the algebra $\mathcal{A}_{FM}$. It is clear that in general one
needs a theory, parallel to cyclic cohomology theory for Hopf
algebras, that works for all bialgebroids endowed with a suitable
additional structure related to an antipode. The goal of this
paper is  to isolate an appropriate  class of bialgebroids, called
here para-Hopf algebroids,
 for which one can
define a cocyclic module and hence a cyclic cohomology theory
extending the cyclic cohomology for Hopf algebras (Definition 2.1
and Theorem 3.1).

Among the properties of Hopf algebras  that helped Connes and
Moscovici in \cite{achm99} (cf. also \cite{achm00} for a survey) to prove that  the structure
discovered  in \cite{achm98} is actually a cocyclic module for any
Hopf algebra endowed with a  modular pair in involution are the
following: the antipode $S$ is an anti-coalgebra map; the twisted
antipode $\tilde{S}$ is a twisted anti-coalgebra map;
$\tilde{S}^2=1$. The first two  properties unfortunately  are not
well defined in the world of bialgebroids. Thus, our first task is
to find good,  i.e. well defined for bialgebroids, equivalents
of these properties.

The main challenge is to find necessary and sufficient conditions, in terms of a single operator called
here a {\it para-antipode}, for the  cyclic operator of the
Connes-Mosocovici module to be well defined and to form a cocyclic module.
One first realizes that the complex defined in \cite{achm01}  is
defined for all bialgebroids endowed with an antialgebra map  $T$
satisfying conditions (PH1) and (PH2) of Definition \ref{one}, and in fact it is always a
cosimplicial module.
After many attempts, we realized that if we just assume that the
 second and
third powers
 of the
cyclic operator satisfy
$$\tau_1^2(h)=h, \quad \text{and} \quad
 \tau_2^3(1_H\ten h)=1_H\ten h, $$
 for all $h \in H$,
 then we have a cocyclic module. These last two conditions are obviously
 necessary as well. We find it remarkable that the third power of the
 cyclic operator comes into the picture and, together with its second
 power, gives necessary and sufficient conditions to have a cocyclic
 module. Note that for Hopf algebras it is the second power of the cyclic operator that gives
 a necessary and sufficient condition to have a cocyclic module \cite{achm99}.

 Our main theorem,  coupled with the fact   that the Connes-Moscovici algebra
 $\mathcal{H}_{FM}$ admits a cocyclic module, implies that $\mathcal{H}_{FM}$
is a para-Hopf algebroid. We also
 provide  a few other examples of  para-Hopf algebroids, including  the
  algebra $A_{\theta}$ of noncommutative
    torus, and the bialgebroid of Example 2.4 below  defined by Connes and Moscovici in
    their study of Rankin-Cohen brackets \cite{cmrc}.

In an earlier version of this paper, of which  a very brief sketch
appeared in our survey article \cite{kr}, a different definition was proposed for a class of
bialgebroids that admit a cocyclic module. The axioms for this structure, called an
extended Hopf algebra there,
are unfortunately very difficult to verify. In particular our claim  there that
the Connes-Moscovici algebra is an extended Hopf algebra seems to be wrong or at least our proof
is not adequate. Our present notion of para-Hopf algebroid on the other
hand seems to be more workable and useful.

    We would like to thank Alain Connes and Henri Moscovici for their   interest and for 
    valuable comments
    and suggestions that played a crucial role in
    the development of our ideas
    specially with regard to the Connes-Moscovici algebra
    $\mathcal{H}_{FM}$.  It is a pleasure to thank Daniel Sternheimer for several valuable suggestions
    that improved our exposition and style. We are also much obliged to a referee who suggested
    that we look into  Example 2.4, as well as 
    suggesting  a new title, and a new name, para-Hopf algebroids,  for the main object
    of study in this paper.

\section{Para-Hopf algebroids}
In this section we first recall the definitions  of  bialgebroids and
Hopf algebroids
  from \cite{lu,Xu}. We
 then define a {\it para-Hopf algebroid} as a bialgebroid endowed with an extra structure that we call a
 {\it para-antipode}.
  Finally we show that several classes of bialgebroids are in fact para-Hopf algebroids. They
  include the Connes-Moscovici algebra
 $\mathcal{H}_{FM}$,  the groupoid algebra of a groupoid with a
 finite number of objects,
 and the algebraic
 noncommutative torus $A_{\theta}$.

Hopf algebroids can be regarded as not necessarily commutative or
cocommutative algebraic analogues of groupoids. Commutative Hopf
algebroids were introduced in \cite{ra}
 as cogroupoid
objects in the category of commutative algebras. Hence the total algebra and the base algebra
are both commutative. Equivalently they can be defined as representable
functors from the category of commutative algebras to the category of
groupoids.
 The main examples of commutative Hopf algebroids are algebras of functions on an algebraic groupoid.

 In \cite{mal} Hopf algebroids were defined where  the total algebra need not be commutative but the base algebra
 is still commutative and the source and target maps    land in the center of
 the total algebra. In the following we first recall a general definition of a bialgebroid
 and a Hopf algebroid, due to Lu, from \cite{lu, Xu} and then define our para-Hopf algebroids. Note that  a general definition of a
 bialgebroid, equivalent to Lu's definition,
  was also given
 by  Takeuchi in \cite{tak}.

Let $k$ be a field of characteristic zero. A {\it bialgebroid}
$(H,R,\alpha,\beta,\Delta,\epsilon)$ over $k$ consists of the
following data:
\begin{itemize}
\item[BA1:] A $k$-algebra $H$, a $k$-algebra $R$, an algebra
homomorphism $\alpha :R \rightarrow H$, and  an anti algebra
homomorphism $\beta : R \rightarrow H$ such that the images of
$\alpha $ and $\beta $ commute in $H$,  i.e. for all $a$, $b$ in
$R$
$$\alpha(a)\beta(b)=\beta(b)\alpha(a). $$
 It follows that $H$ has an  $R$-bimodule structure defined by
\begin{center}
$axb=\alpha(a)\beta(b)x$~~~~ $\forall a,b\in R$, ~$x\in H.$\\
\end{center}
In particular the bimodule tensor product $H\otimes_R H$ is defined and
is an $(R, R)$-bimodule. Similarly for $H\otimes_R H\otimes_R H$ and
higher bimodule tensor products.
$H$ is called the {\it total algebra}, $R$ the {\it base algebra},
$\alpha$ the {\it source map},
and $\beta$
the {\it target map} of the bialgebroid.
\item[BA2:]A coproduct, i.e.  an $(R,R)$-bimodule map $\Delta : H\rightarrow H\otimes_RH$
which  satisfies the following conditions:
\begin{itemize}
\item[cp1)]$\Delta(1)=1\otimes_R1$
\item[cp2)]Coassociativity :
$$(\Delta \otimes_R id_H)\Delta = (id_H\otimes_R\Delta)\Delta :H \rightarrow H\otimes_RH\otimes_RH,$$
\item[cp3)] Compatibility with the  product: for all $a$, $b\in H$ and $r\in R$,
$$\Delta(a)(\beta(r)\otimes 1-1\otimes\alpha(r) )=0 ~~~ \text{in}~ H\otimes_RH, $$
$$\Delta (ab)=\Delta(a)\Delta(b).$$
In the first relation the natural right action of $H\otimes H$ on
$H\otimes_R H$ defined by $(a\otimes_R b)(a'\otimes
b')=aa'\otimes_R bb'$ is used. While $H\otimes_R H$ need not be an
algebra, it can be easily checked that the left annihilator of the
image of $\beta \otimes 1-1\otimes \alpha $ is an algebra. Hence,
by the first relation,  the multiplicative property of $\Delta$
makes sense.
\end{itemize}
\item[BA3:] A counit, i.e. an $(R,R)$-bimodule  map $\epsilon:H\rightarrow R$ satisfying,
\begin{itemize}
\item[cu1)]
$\epsilon(1_H)=1_R$
\item[cu2)] $(\epsilon \otimes_Rid_H )\Delta =(id_H\otimes_R\epsilon)\Delta =id_H:H\rightarrow  H$
\end{itemize}
\end{itemize}

In this paper we suppress the summation notation and write
$\Delta(h)=h^{(1)}\otimes_R h^{(2)}$ to denote the coproduct of a
bialgebroid. Similarly,  $\Delta^n(h)=h^{(1)}\otimes_R\cdots
\otimes_Rh^{(n+1)}$ denotes the higher iterations of the coproduct.

 Although it is  not used in this paper, for reader's convenience in comparing the definitions,
 we recall that
 a  bialgebroid $(H,R,\alpha , \beta, \Delta ,\varepsilon)$
is called  a {\rm Hopf algebroid}
 if there is a bijective map $S:H\rightarrow H$, called antipode, which is an antialgebra map
 satisfying the following conditions:
\begin{itemize}
\item[1)]\label{one}$S\beta = \alpha $.
\item[2)]\label{three}$m_H(S\otimes id)\Delta =\beta \epsilon
S:H\rightarrow H$,
where $m_H : H\otimes H \rightarrow H$ is the multiplication map of
$H$.
\item[3)]\label{three} There exists a linear map $\gamma:H\otimes_RH \rightarrow H\otimes H$ satisfying\\
$\pi\circ\gamma =id_{H\otimes_RH} :H\otimes_RH\rightarrow H\otimes_RH$ and $m_H(id\otimes S)\gamma
\Delta =\alpha \epsilon :H\rightarrow H$\\
where $\pi:H \otimes H \rightarrow H \otimes _RH $ is the  natural projection.
\end{itemize}

Note that while the operator $m_H(S\otimes id)\Delta :H
\rightarrow H$ is well defined, i.e.,
 is independent of the choice of any section for the projection map
 $H\otimes_R H \rightarrow H\otimes H$,
the operator $m_H(id\otimes S) \Delta $ is not well defined and
one has to fix a linear section $\gamma $ for $\pi$ first. Also, unlike
Hopf algebras, the antialgebra property of $S$ does not follow from axioms $HA1)-HA3)$,
and has to be assumed. Since the operator $H\otimes_R H \rightarrow
H\otimes_R H$, $x\otimes_Ry\mapsto y\otimes_Rx$ is not even well
defined, the anti-coalgebra property of $S$ does not even make sense.

\begin{definition}
 A bialgebroid $(H,R,\alpha , \beta, \Delta ,\epsilon)$
is called  a {\rm Para-Hopf algebroid}
 if there is an antialgebra map $T:H\rightarrow H$, called a {\rm para-antipode},
 satisfying the following conditions:
\begin{itemize}
\item[PH1)]\label{one}$T\beta = \alpha $.
\item[PH2)]\label{three}$m_H(T\otimes id)\Delta =\beta \epsilon
T:H\rightarrow H$, where $m_H : H\otimes H \rightarrow H$ is the
multiplication map of $H$.
 \item[PH3)] $T^2=id_H$, and for all $h
\in H$
\begin{equation}\label{anticoalgebra}
 T(h\sw{1})\sw{1}h\sw{2}\ten
T(h\sw{1})\sw{2}=1\otimes_RT(h).
\end{equation}
\end{itemize}
\end{definition}

\begin{remark}
It follows from Theorem 3.1 that in  terms of the cyclic operator $\tau$  in Theorem 3.1, axiom $PH3$
 can be expressed as:
$$\tau_1^2(h)=h, \quad \text{and} \quad
 \tau_2^3(1_H\ten h)=1_H\ten h, $$
for all $h\in H$.
\end{remark}

We label the result of the following lemma (cu3) because some authors assume it as an additional
axiom  for the
 counit in the definition of a bialgebroid.  We  need (cu3) in several
 proofs in this paper,
 especially in the proof of Theorem \ref{main}.
\begin{lem}
Let $(H,R)$ be  a Para-Hopf algebroid. Then  for all $h,g$ in $H$ we have

${\rm (cu3)}\quad \epsilon(hg)=\epsilon(h\alpha(\epsilon(g)))=\epsilon(h\beta(\epsilon(g))).$
 \end{lem}
 \begin{proof}

 By using  (PH1), (PH2)  and $T^2=Id_H$, we have:
 \begin{align*}
 &\epsilon(hg)=\epsilon(T^2(hg))=T(T(hg)\sw{1})T(hg)\sw{2}=T(T(h)\sw{1})T(T(g)\sw{1})T(g)\sw{2}T(h)\sw{2}\\
 &=T(T(h)\sw{1})\beta(\epsilon(g)) T(h)\sw{2}=T(T(h\beta(\epsilon(g)))\sw{1}) T(h\beta(\epsilon(g)))\sw{2}=\epsilon(h\beta(\epsilon(g))).
  \end{align*}
  The other equality can be proven the same way.
  \end{proof}

Throughout this paper we make use  of  the following  two
actions. First, the right action of $H^{\ot n}$ on $H^{\ten n}$
defined by
\begin{equation}
(h_1\ten \dots \ten h_n)\cdot(g_1\ot \dots \ot g_n)=(h_1g_1\ten
 h_2g_2\ten\dots \ten h_ng_n).
\end{equation}
It is evident that this  action is well defined  because the  $R$-bimodule structure
of $H$ is defined by using  $\alpha$, $\beta$ and  left
multiplication.  The next action is  the left action of $H$ on
$H^{\ten n}$ defined by
\begin{equation}\label{diagonal}
h\vartriangleright(g_1\ten \dots \ten g_n)= h\sw{1}g_1\ten\dots
\ten h\sw{n}g_n.
\end{equation}
This is a well defined  action because of  property $cp3)$ of the
coproduct of a bialgebroid. One can show that $H^{\ten n}$ is an $H-H^{\ot n}$
bimodule. We use this bimodule structure in some proofs in this
paper.

 The following lemma will prove useful in verifying that certain
examples satisfy  condition (1) of definition 2.2.
\begin{lem} \label{multi}Condition (1) above is multiplicative. That is, if it is
satisfied  by $h$ and $g$ then it is satisfied by $hg$.
\end{lem}
\begin{proof}
Let  $h$ and $g$ satisfy  (\ref{anticoalgebra}).
 We have  \\ $ T((hg)\sw{1})\sw{1}(hg)\sw{2}\ten
T((hg)\sw{1})\sw{2}$\\$=T((hg)\sw{1})\vartriangleright ((hg)\sw{2}\ten
1_H)$\\
$=T(h\sw{1}g\sw{1})\vartriangleright(h\sw{2}g\sw{2}\ten
1_H)$\\$=(T(g\sw{1})T(h\sw{1}))\vartriangleright(h\sw{2}g\sw{2}\ten
1_H)$\\
$=T(g\sw{1})\vartriangleright(T(h\sw{1})\vartriangleright(h\sw{2}g\sw{2}\ten
1_H))$\\$
=T(g\sw{1})\vartriangleright(T(h\sw{1})\sw{1}h\sw{2}g\sw{2}\ten
T(h\sw{1})\sw{2})$\\
$=T(g\sw{1})\vartriangleright ((T(h\sw{1})\sw{1}h\sw{2}\ten
T(h\sw{1})\sw{2})\cdot(g\sw{2}\ot 1_H))$\\$=T(g\sw{1})\vartriangleright((1_H\ten T(h))\cdot (g\sw{2}\ten
1_H))$\\
$=T(g\sw{1})\vartriangleright( (g\sw{2}\ten 1_H)\cdot(1_H\ot
T(h)))$\\
$=(T(g\sw{1})\vartriangleright( g\sw{2}\ten 1_H))\cdot(1_H\ot
T(h))$\\$=(1_H\ten T(g))\cdot(1_H\ot T(h))$\\
$=1_H\ten T(g)T(h)=1_H\ten T(hg).$\\
This shows that $hg$ satisfies (1) as well.
\end{proof}

The following proposition shows that for Hopf algebras our para-antipodes
are simply twisted antipodes in the sense of Connes and
Moscovici \cite{achm98}. We would like to emphasize that no analogue of
this result exists in the world of bialgebroids.

\begin{prop} Let $\mathcal{H}$ be a  Hopf algebra over $k$ with   antipode $S$.
 Then   $T : \mathcal{H}\rightarrow \mathcal{H}$ is an antialgebra map  and satisfies
  (\ref{anticoalgebra})
   if and only if $T=\delta\ast S$, where
$\ast$ denotes  the convolution multiplication and
$\delta:\mathcal{H}\rightarrow k$ is an algebra map.
 \end{prop}
 \begin{proof}
 Let  $\delta$ be  an algebra map.  It is obvious that $T=\delta\ast S: \mathcal{H}\rightarrow
 \mathcal{H}$ is an antialgebra map. We have
 $$T(h\sw{1})\sw{1}\ot T(h\sw{1})\sw{2}=\delta(h\sw{1})S(h\sw{2})\sw{1}h\sw{3}\ot
 S(h\sw{2})\sw{2}$$
 $$=\delta(h\sw{1})S(h\sw{3})h\sw{4}\ot
 S(h\sw{2})=1_\Hc \ot \delta(h\sw{1})S(h\sw{2})=1_\Hc\ot T(h).$$

 On the other hand, let $T$ be  an antialgebra map that  satisfies
 condition (\ref{anticoalgebra}).  We define  $\delta=\epsilon \circ T :\Hc\rightarrow
 k$.  It is evident that $\delta $ is an algebra map. We  verify that $\delta\ast
 S=T$. Indeed, by our assumption  we have\\
 $$T(h\sw{1})\sw{1}h\sw{2}\otimes
T(h\sw{1})\sw{2}=1\otimes T(h).$$
 Applying $m\circ (S\ot id_\Hc)$, where $m$ denotes the
 multiplication map of $\Hc$, to both  sides  of the above
 equation, we obtain
 $$S(h\sw{2})S(T(h\sw{1})\sw{1})T(h\sw{1})\sw{2}=T(h),$$
or, equivalently,  $\delta(h\sw{1})S(h\sw{2})=T(h).$
  \end{proof}

  We give a few examples of para-Hopf algebroids.
\begin{example}\label{reza}{\rm
Let $\mathcal{H}$ be a  Hopf algebra over a field $k$,
$\delta:\mathcal{H}\rightarrow k$ an algebra map, and $\widetilde{S}_{\delta}=\delta\ast S$
 the $\delta$-twisted antipode defined by $\widetilde{S}_{\delta}(h)=\delta(h^{(1)})S(h^{(2)}).$
  Assume that $\widetilde{S}_{\delta}^2=id$.
  Let $R$ be any algebra over $k$. We define a para-Hopf algebroid as
  follows.
  Let
  $H=R\otimes \mathcal{H}\otimes R^{op}$, where $R^{op}$ denotes the opposite algebra of $R$.
 One can check that with the following structure $(H,R)$ is an extended
 Hopf algebra:

\begin{eqnarray*}
\alpha(a)&=&a\otimes 1\otimes1 \\
\beta(a)&=&1\otimes 1 \otimes a\\
\Delta(a\otimes h\otimes b)&=& a\otimes h^{(1)}\otimes 1\otimes_R 1\otimes h^{(2)}\otimes b\\
\epsilon(a\otimes h\otimes b)&=&\epsilon(h)ab\\
T(a\otimes h\otimes b)&=&(b\otimes \widetilde{S}_\delta(h)\otimes
a).
\end{eqnarray*}}
\end{example}

\begin{example}\label{reza3}{\rm
Let $\mathcal{G}$ be a groupoid over a finite base. Equivalently, $\mathcal{G}$
is  a category with a finite set of  objects, such that
 each morphism is invertible. The groupoid algebra of $\mathcal{G}$, denoted by $H=k\mathcal{G}$,
  is freely generated over $k$ by
morphisms $g\in \mathcal{G}$ with unit $1=\sum_{X\in \mathcal{O}bj(\mathcal{G}) }id_X$. The
 product  of two
morphisms is equal to their composition  if the latter is defined and $0$ otherwise. We show
that $k\mathcal{G}$
is an
extended  Hopf algebra over the base algebra
$R=k\mathcal{S}$, where $\mathcal{S}$ is the subgroupoid of $\mathcal{G}$
whose objects  are those of
 $\mathcal{G}$ and $\mathcal{M}or(X,Y)=id_X$ whenever $X=Y$ and $\emptyset$ otherwise.
 The relevant  maps are defined as follows:
 $\alpha=\beta: R\hookrightarrow H$ is the
 natural embedding, and
  $$\Delta(g)=g\otimes_Rg, \   \epsilon(g)=id_{target(g)}, \  T(g)=g^{-1},$$
  for any
  $g\in \mathcal{G}$.
 To show that it is an para-Hopf algebroid, we see that  all conditions are obvious except
 possibly the condition
  (\ref{anticoalgebra}). To check this, we compute $\tau_2^3(1_H\ten g)$.
  We have $\tau_2(1_H\ten g)=g\ten 1_H$, and  $\tau_2(g\ten 1_H)=g^{-1}\ten g^{-1}$. Hence
  $\tau_2^3(1_H\ten g)=gg^{-1}\ten g=1_H\ten g$. }
\end{example}

\begin{example}{\rm
 It is known that the algebraic quantum torus $A_{\theta}$  is not
 a Hopf algebra
 although it is a deformation
 of the  Hopf algebra of Laurent polynomials in two variables. We show that $A_{\theta}$ is a
 para-Hopf algebroid. Recall  that  $A_{\theta} $ is
 the unital $\mathbb{C}$-algebra generated by two invertible elements $U$ and $V$ subject to the
   relation $UV=qVU$, where
 $q=e^{2\pi i{\theta}}$ and $\theta \in \mathbb{R}$. Let $R=\mathbb{C}\mathbb\lbrack U, U^{-1}\rbrack$
 be the algebra of Laurent
 polynomials embedded in $A_{\theta}$. Let $\alpha = \beta : R\rightarrow A_{\theta}$ be the
 natural embedding. Define the coproduct
 $\Delta : A_{\theta}\rightarrow A_{\theta} \otimes_R A_{\theta} $ by
  $$\Delta(U^nV^m)=U^nV^m\otimes_R V^m $$
   and the counit $\varepsilon : A_{\theta} \rightarrow R $ by
   $\varepsilon(U^nV^m)=U^n$. Define the extended antipode $T:A_{\theta} \rightarrow A_{\theta}$
   by
   $$T(U^nV^m)=q^{nm}U^nV^{-m}.$$ We  check   that
   $A_{\theta}$ is
   a para-Hopf algebroid over $R$.
   Among the axioms we just check the validity of the condition (\ref{anticoalgebra}).
   By   Lemma \ref{multi}, it is enough  to  check this condition just for the
   generators $U$ and $V$.
   We have
   $$\tau_2^3(1_h\ten U)=\tau^2_2(U\ten 1_H)=\tau_2(U\ten 1_H)=(U\ten
   1_H)=1_H\ten U,$$
   and
   $$\tau_2^3(1\ten V)=\tau_2^2(V\ten 1_H)=\tau(V^{-1}\ten V^{-1})=1_H\ten V.$$
    }
\end{example}

Next we show that   the Connes-Moscovici
algebra $\mathcal{H}_{FM}$, introduced in \cite{achm01},
 is a para-Hopf algebroid. In fact it is already shown in \cite{achm01} that $\mathcal{H}_{FM}$
is a bialgebroid and an antialgebra map  $\tilde{S}: \mathcal{H}_{FM} \rightarrow \mathcal{H}_{FM}$
is defined
such that $\tilde{S}^2=id$,  $\tilde{S}\beta = \alpha$, and $m_H(\tilde{S}\otimes id)\Delta =
\beta \epsilon
\tilde{S}:\mathcal{H}_{FM}\rightarrow \mathcal{H}_{FM}$. All we have to do then
is to check that condition (1) is satisfied for $T=\tilde{S}$.
 First let us briefly recall the
definition of
$\mathcal{H}_{FM}$ from \cite{achm01}.

 Let $M$ be an smooth n-dimensional manifold that admits a finite open
cover by coordinate charts and let $FM$ denote the frame bundle of
$M$. A local diffeomorphism of $FM$ is called a  prolonged
diffeomorphism if it is the natural prolongation of a local
diffeomorphism of $M$. Let $\mathcal{G}$ denote the set of germs
of prolonged local diffeomorphisms of $FM$. Then $\mathcal{G}$ is
a smooth \'etale groupoid with $FM$ as its set of objects. Let
$\mathcal{A}=\mathcal{A}_{FM}:= C^{\infty}_c(\mathcal{G})$ denote
the smooth convolution  algebra of $\mathcal{G}$. Its elements are
linear combinations of elements of the form  $fU_{\phi}^*$, with
$f\in C^{\infty}_c(Dom \tilde{\phi})$. Here  $\tilde{\phi}$
denotes the prolongation of  a local diffeomorphism $\phi $ of $M$
and the asterisk  denotes the inverse.  The product is defined by
$$ f_1 U_{\phi_1}^* \cdot f_2 U_{\phi_2}^*= f_1 (f_2\circ \tilde{\phi}_1) U_{\phi_2 \phi_1}^*. $$
Let $\mathcal{R}=\mathcal{R}_{FM}: =C^{\infty}(FM)$ denote  the
algebra of smooth functions on the frame bundle $FM$.
$\mathcal{R}$ acts on $\mathcal{A}$ by left and right
multiplication operators: $\alpha (b) (fU_{\phi}^*)= b \cdot f
U_{\phi}^*$ and $\beta (b) (fU_{\phi}^*)=fU_{\phi}^* \cdot b =b
\circ \phi \cdot fU_{\phi}^*$, for all $b\in \mathcal{R}.$ It is
easily seen that $\alpha : \mathcal{R}\rightarrow End
(\mathcal{A})$ is  an algebra map, $\beta : \mathcal{R}\rightarrow
End (\mathcal{A})$ is an antialgebra map, and the images of
$\alpha$ and $\beta$ commute. We also have the action of  vector
fields on $FM$ on the algebra
 $\mathcal{A}$ by the formula $Z (fU_{\phi}^*)= Z(f)
U_{\phi}^*$, where $Z$ is a vector field. Note that while  vector
fields act by derivations
 on functions on the frame bundle, their action on $\mathcal{A}$ does not satisfy the
 derivation
 property. In fact the failure of derivation property is responsible for   the
 non-cocommutativity of the coproduct
 of $\mathcal{H}_{FM}$.

Let $\mathcal{H}_{FM}\subset End (\mathcal{A})$ be the subalgebra of the
algebra of
linear operators on $\mathcal{A}$ generated by the images of
$\alpha$, $\beta$ and actions of vector fields as above. Its elements
are called {\it transverse differential operators} on the \'etale groupoid
$\mathcal{G}$.

 It is shown  in \cite{achm01} that $\mathcal{H}_{FM}$
   is a free $\mathcal{R}\otimes \mathcal{R}$-module.
    In fact  fixing a torsion free connection on $FM$, one obtains a
    Poincar\'e-Birkhoff-Witt-type basis for $\mathcal{H}_{FM}$ over $\mathcal{R}\otimes \mathcal{R}$
    as follows. Let $\{Y^i_{j}\}$ denote the fundamental vertical vector
    fields corresponding to the standard basis of $gl(n, \mathbb{R})$ and $X_1, \cdots X_n$ denote the
    standard horizontal vector fields corresponding to the standard basis of $\mathbb{R}^n$. These $n^2 +n$ vector
    fields form a basis for the tangent space of $FM$ at all points. It is shown in \cite{achm01} that
    the operators
    $Z_I\cdot\delta_{\kappa}$, where
    $$Z_I=X_{i_1}\cdots X_{i_p}Y^{j_1}_{k_1}\cdots Y^{j_q}_{k_q},
    \ \ \delta_{\kappa} =\delta^{i_1}_{j_1k_1; \ell^{1}_1 \cdots \ell^{1}_{p_1}}\cdots
    \delta^{i_r}_{j_rk_r; \ell^{r}_1 \cdots \ell^{r}_{p_r} },$$
    $$\delta^{i}_{jk; \ell_1 \cdots \ell_p}=[X_{\ell_r}\cdots [X_{\ell_1},
    \delta^i_{jk}]\cdots],$$
    form a basis for $\mathcal{H}_{FM}$  over $\mathcal{R}\otimes
    \mathcal{R}$.
(See Proposition 3 and Lemma 2 in \cite{achm01} for  precise range of multi-indices $I$ and $\kappa$ as well
as the definition of the
operators $\delta^i_{jk}$.)

     A coproduct $\Delta$ and an antialgebra map $\widetilde{S}$ with $\widetilde{S}^2=id$ are already
    defined in \cite{achm01} and all the identities of a bialgebroid as well as axioms
    $PH1, PH2$ in Definition 2.2 are verified.
    We  check that condition (\ref{anticoalgebra}) is satisfied as well.

    \begin{lem} With $T=\tilde{S}$ condition (\ref{anticoalgebra}) is satisfied for the
    Connes-Moscovici algebra $\mathcal{H}_{FM}$.
    \end{lem}
    \begin{proof}
    Thanks to Lemma {\ref{multi}}, we just need to check the
    condition (\ref{anticoalgebra}) for the generators of
    $H=\mathcal{H}_{FM}$. Let $R=\mathcal{R}$. We just check the validity of this condition
    for generators
    $X_k$, the rest being   straightforward to check.
     We know  that $$\Delta(X_k)=X_k\ten 1_H +1_H\ten X_k +\delta^i_{jk}\ten
     Y^j_i.$$ So the right hand side of the condition
     (\ref{anticoalgebra}) converts to
     $$\ts(X_k)\sw{1}\ten \ts(X_k)\sw{2}+ X_k\ten 1_H +\ts (\delta^i_{jk})\sw{1}
     Y^j_i\ten\ts(\delta^i_{jk}).$$
    Now by replacing  $\ts(X_k)$ and $\ts(\delta^i_{jk})$ in the above expression by their
    equals    $-X_k+\delta^i_{jk}Y^j_i$ and
    $-\delta^i_{jk}$   respectively, and using
     the fact that $\Delta $ is multiplicative,  we    find that the above
     expression is equal to
    \begin{align*}
    &X_k\ten 1_H-1_H\ten X_k-\delta^i_{jk}\ten Y^j_i+\delta^i_{jk}Y^j_i\ten 1_H+\delta^i_{jk}\ten Y^j_i\\
    & + Y^j_i\ten \delta^i_{jk}
     +1_H\ten \delta^i_{jk}Y_i^j + X_k\ten 1_H-\delta^i_{jk}Y^j_i\ten 1_H-Y^j_i\ten
     \delta^i_{jk}.
     \end{align*}
     After cancelling  the identical  terms with
     opposite signs  we obtain
    $-1_H\ten X_k+1_H\ten \delta^i_{jk}Y_i^j$, which is  $1_H\ten \ts(X_k).$
    \end{proof}

\begin{remark}
The fact that the Connes-Moscovici algebra $\mathcal{H}_{FM}$ is a
para-Hopf algebroid can also be directly derived by combining
Theorem 3.1 with the fact, proved in \cite{achm01}, that
$\mathcal{H}_{FM}$ affords a cocyclic module.
\end{remark}

If $M=\mathbb{R}^n$ is the flat Euclidean space,
then $ \mathcal{H}_{FM} =\mathcal{R}\otimes \mathcal{H}_n \otimes
\mathcal{R}$, where $\mathcal{H}_n$ is the Connes-Moscovici Hopf
algebra in dimension $n$. The
  para-Hopf algebroid structure on $\mathcal{H}_{FM}$ is  induced from the Hopf algebra structure
  on $\mathcal{H}_n$ as in Example \ref{reza}.
\begin{example}
\rm{Let $\Hc$ be a Hopf  algebra and $\cp$ be  a left
$\mathcal{H}$-module algebra. We generalize  Example \ref{reza} by
turning the  double crossed product algebra $H:=\mathcal{P}\rtimes
\mathcal{H}\ltimes \mathcal{P}^{\text{op}}$, introduced in
\cite{cmrc}, into a Para-Hopf algebroid. To this end let us first
recall from \cite{cmrc}
 its   algebra structure.  Equipped with the
 following multiplication and  $1\rtimes 1\ltimes 1$ as its unit,    $ \cp\ot\Hc\ot \cp$ is a
 unital  associative algebra:
 \begin{equation*}
 (P_1\rtimes h_1\ltimes Q_1)\cdot (P_2\rtimes h_2\ltimes Q_2)=P_1h_1\sw{1}(P_2)\rtimes h_1\sw{2}h_2\ltimes h_1\sw{3}(Q_2)Q_1.
 \end{equation*}
 The source and target maps   are defined by  $\alpha: \cp\rightarrow H$ defined by $\alpha(P)=P\rtimes 1\ltimes 1  $ and $\beta:\cp
  \rightarrow
  H$ defined by $\beta(Q)=1\rtimes 1\ltimes Q$.
  The comultiplication  $\Delta:H\rightarrow H\ot_\cp H $ is defined by
  $$\Delta(P\rtimes h\ltimes Q)=P\rtimes h\sw{1}\ltimes 1\ot_\cp 1\rtimes
   h\sw{2}\ltimes Q.$$
   The above data together  with the counit $\epsilon:\cp \rightarrow H$ defined  by $\epsilon(P\rtimes h\ltimes Q)
 =P\epsilon(h)Q$
turn  $H$ into a bialgebroid. Now if $\Hc$ admits a character such that
$\widetilde{S_\delta}^2=id_\Hc$ then  we  define   a para-antipode for $H$ and show
 that $H$ is a para-Hopf algebroid over $\cp$. Let $T: H\rightarrow H$ be defined by
 \begin{equation*}
 T(P\rtimes h\ltimes Q)=S(h\sw{3})(Q)\rtimes S(h\sw{2})\ltimes \widetilde{S_\delta}(h\sw{1})(P)
 \end{equation*}
 Evidently $T^2=id$, and $T$ is an antialgebra map.   The other  identities are straightforward to check and we leave them to the reader except the
   crucial identity (\ref{anticoalgebra}). Let $\mathfrak{h}=P\rtimes h\ltimes Q\in H.$ we have
\begin{eqnarray*}
   T(\mathfrak{h}\sw{1})\sw{1}\mathfrak{h}\sw{2}\ot_\cp
   T(\mathfrak{h}\sw{1})\sw{2}
   &=&(1\rtimes S(h\sw{3})\ltimes 1)\cdot(1\rtimes h\sw{4}\ltimes
   Q)\ot_\cp \\
   & &1\rtimes S(h\sw{2})\ltimes
  \widetilde{S_\delta}(h\sw{1})(P)\\
   &=&1\rtimes 1\ltimes S(h\sw{3})(Q)\ot_\cp 1\rtimes S(h\sw{2})\ltimes
   \widetilde{S_\delta}(h\sw{1})(P) \\
    &=&1\rtimes 1\ltimes 1 \ot_\cp S(h\sw{3})(Q)\rtimes S(h\sw{2})\ltimes \widetilde{S_\delta}(h\sw{1})(P)\\
    &=&1_H\ot_\cp T(\mathfrak{h}).
  \end{eqnarray*}}
\end{example}

\section{Cyclic cohomology of para-Hopf algebroids}
 In \cite{achm01},   Connes and  Moscovici used the natural action
  of $\mathcal{H}_{FM}$ on the algebra $ \mathcal{A}_{FM}$ and an invariant
  faithful trace $Tr$ on $ \mathcal{A}_{FM}$ to define  a cocyclic
  module for $\mathcal{H}_{FM}$. More precisely, they showed that the
  maps
  $$\gamma_{Tr}:\mathcal{H}_{FM}^{\otimes_R(n+1)} \longrightarrow
  Hom_{\mathbb{C}} (\mathcal{A}_{FM}, \mathbb{C}),$$
  defined by
  $$\gamma_{Tr} (h_1\otimes_R \cdots \otimes_R h_n)(a_0\otimes
  \cdots \otimes a_n)=Tr (a_0 h_1(a_1)\cdots h_n(a_n))$$
  are   linear isomorphisms for each $n\geq 0$ and their images form  a cocyclic submodule of the cocyclic module of the
  algebra $\mathcal{A}_{FM}$. This cocyclic submodule was then transferred,
  via $\gamma_{Tr}$, to a cocyclic module based on $\mathcal{H}_{FM}$.

In this section our aim is to show that the formulas discovered by Connes and Moscovici define a cocyclic module for
any para-Hopf algebroid and provide several computations.

 Let $(H,R)$ be a bialgebroid. It is easily checked that
 the following module $H_\natural$ is a cosimplicial module.   We put
  $$H_\natural^0=R, \text{ and} \;H_\natural^n=H\otimes_R H\otimes_R
  \dots \otimes_RH \qquad (n\text{-fold tensor product}).$$
   The cofaces $\delta_i$ and  codegeneracies  $\sigma_i$  are defined by:
\begin{eqnarray*}
\delta_0(a)=\alpha(a),~ \delta_1(a)=\beta(a)&& \text{for all}~ a\in R=H^0_\natural\\
\delta_0(h_1\otimes_R\dots \otimes_Rh_n)&=&1_H\otimes_R h_1\otimes_R\dots \otimes_Rh_n \\
\delta_i(h_1\otimes_R\dots \otimes_Rh_n)&=& h_1\otimes_R\dots\otimes_R\Delta (h_i)\otimes_R\dots \otimes_Rh_n
  \;\;\text{for}\;\;1\leq i\leq n \\
\delta_{n+1}(h_1\otimes_R\dots \otimes_Rh_n)&=&h_1\otimes_R\dots \otimes_Rh_m\otimes_R 1_H \\
\sigma_i(h_1\otimes_R\dots \otimes_Rh_n)&=& h_1\otimes_R\dots\otimes_R\epsilon(h_{i+1})\otimes_R\dots \otimes_Rh_n
 \;\;\text{for}\;\;0\leq i\leq n. \\
\end{eqnarray*}
Let $T: H \rightarrow H$ be  an antialgebra map such that axioms $PH1$  and
$PH2$ of Definition 2.2 are satisfied. That is
  $T\beta =\alpha$
and $m_H(T\otimes id)\Delta =\beta \epsilon
T:H\rightarrow H$. Then one can check that
the cyclic operator $\tau$ of Connes and Moscovici defined by
$$\tau_n(h_1\otimes_R\dots \otimes_Rh_n
)=\Delta^{n-1}T(h_1)\cdot(h_2\otimes\dots \otimes
h_n\otimes 1_H),$$
is well defined. Note that the
 right action of $H^{\otimes n}$ on $H^{\otimes_Rn}$ by right
  multiplication is used.
    The question of when $H_\natural$, endowed with the above cyclic operator, is a cocyclic module is
    completely answered by the following theorem.
\begin{theorem}\label{main}

Let $(H,R)$ be a bialgebroid and $T:H \rightarrow H$ be as above. Then
$H_\natural$ is a cocyclic module if and only if $(H, R, T)$ is a para-Hopf algebroid, that is
 $T^2=id_H$, and for all $h\in H,$
 $$T(h\sw{1})\sw{1}h\sw{2}\ten T(h\sw{1})\sw{2}=1\otimes_RT(h).$$
\end{theorem}
\begin{proof}

We should verify the following identities: \\
$\delta_j\delta_i = \delta_i\delta_{j-1} ~~~ i<j$\\
$\sigma_j\sigma_i = \sigma_i\sigma_{j+1} \hspace{15pt} i\le j   $\\
$\sigma_j\delta_i =
\begin{cases}
\delta_i\sigma_{j-1}\;\;\; \text{$i<j$}
  \\
id \;\;\;\text{$i=j$ or $i=j+1$} \\

\delta_{i-1}\sigma_j\;\;\; \text{$i>j+1$}
\end{cases}
$\\
$\tau_{n+1}\delta_i = \delta_{i-1}\tau_n \;\;\; 1\le i\le n$\\
 $\tau_{n-1}\sigma_i = \sigma _{i-1}\tau_n\;\;\; 1\le i\le n$\\
$\tau_n^{n+1} = \mbox{id} _n . $\\

The cosimplicial relations are not hard to prove.
 We just verify the cyclic relations $\tau_n^{n+1} = \mbox{id} _n$ and $\tau\delta_1=\delta_0\tau$,
and leave the others to the reader.

This is evident for $n=1$  because  $\tau_1=T$. Let $n\geq 2 $,  and
define
\begin{center}
$ \Phi_n:H\ten H \longrightarrow   H^{\ten n},$\\
$h\ten g \mapsto \tau_n^2(h\ten g\ten 1\ten\dots\ten 1).$
\end{center}
 We have
$$\tau_n^2(h_1\otimes_R\dots \otimes_Rh_n)=\Phi_n(h_1\ten
h_2)\cdot(h_3\ot\dots\ot h_n\ot 1_H\ot 1_H).$$
 On the other hand we can compute $\Phi_n(h_1\ten
h_2)$ in terms of the diagonal action of  $H$ on $H^{\ten n}$, i.e.,
$$\Phi_n(h_1\ten h_2)=T(h_2)\rhd \Phi_n(h_1\ten 1_H).$$\\
Using condition (\ref{anticoalgebra}), one has $$\Phi_n(h_1\ot
1_H)=1_H\ten 1_H\dots \ten 1_H\ten h_1.$$
We can simplify
$\tau_n^2(h_1\ten \dots\ten h_n)$ as follows:
\begin{align*} &\tau_n^2(h_1\ten\dots\ten h_n)=\\
&T(h_2)\sw{1}h_3\ten \dots\ten T(h_2)\sw{n-2}h_n\ten
T(h_2)\sw{n-1}\ten T(h_2)\sw{n}h_1.
\end{align*}
 By repeating the same argument as above we obtain:
 \begin{align*}
& \tau_n^n(h_1\ten\dots\ten h_n)= T(h_n)\sw{1}\ten
T(h_n)\sw{2}h_1\ten \dots\ten \ten T(h_n)\sw{n}h_{n-1}.
\end{align*}
 Applying $\tau_n$ to both side, we obtain:\\
 $$\tau_n^{n+1}(h_1\ten\dots\ten h_n)=\Phi(h_n\ot 1_H)\cdot( h_1\ot\dots \ot h_{n-1}\ot
 1_H).$$
 Now since $\Phi(h_n\ot
1_H)=1_H\ten 1_H\dots \ten 1_H\ten h_n$, we have
$$\tau_n^{n+1}=id_n.$$

Next,  we check the  identity  between $\tau_n$ and $\delta_1$,
i.e. $\tau_{n+1}\delta_1=\delta_0\tau_n$.
\begin{align*}
&\tau_{n+1}\delta_1(h_1\ten h_2\dots\ten h_n)=\tau_{n+1}(h_1\sw{1}\ten h_2\sw{2}\ten h_2\ten\dots \ten h_n)\\
&=T(h\sw{1})\sw{1}h_1\sw{2}\ten T(h_1\sw{1})\sw{2}h_2\ten \dots
\ten T(h_1\sw{1})\sw{n}h_n\ten T(h_1\sw{1})\sw{n+1}\\
&=(T(h\sw{1})\sw{1}h_1\sw{2}\ten T(h_1\sw{1})\sw{2}\ten \dots
\\
& \dots \ten T(h_1\sw{1})\sw{n}\ten
T(h_1\sw{1})\sw{n+1})\cdot(1_H\ot
h_2\ot\dots\ot h_n\ot 1_H)\\
&=((id_H\ten\Delta\sw{n-2})(T(h\sw{1})\sw{1}h\sw{2}\ten
T(h\sw{1})\sw{2}))\cdot (1_H\ot h_2\ot\dots\ot h_n\ot 1_H)\\
&=((id_H\ten\Delta\sw{n-2})(1_H\ten T(h_1)))\cdot (1_H\ot h_2\ot\dots\ot h_n\ot 1_H)\\
&=1_H\ten T(h_1)\sw{1}h_2\ten\dots\ten T(h_2)\sw{n-1}h_n\ten
T(h)\sw{n}\\
&=\delta_0\tau_n(h_1\ten h_2\dots\ten h_n).
\end{align*}

Finally we check the relation between $\tau_n$ and $\sigma_i$. Using Lemma \ref{cu3} one has:
 \begin{align*}
 &\tau_{n-1}\sigma_i(h_1\ten \dots \ten h_n)=\tau_{n-1}(h_1\ten \dots \ten \beta(\epsilon(h_i))h_{i-1}\ten\dots \ten h_n)\\
 &=S(h_1)\sw{1}h_2\ten\dots\ten S(h_1)\sw{i-2}\beta(\epsilon(h_i))h_{i-1}\ten\dots\ten S(h_1)\sw{n-2}h_n\\
 &=S(h_1)\sw{1}h_2\ten\dots\ten \beta(\epsilon((S(h_1)
 \sw{i-2}\beta(\epsilon(h_i)))\sw{2}))(S(h_1)\sw{i-2}\beta(\epsilon(h_i)))\sw{1}h_{i-1}\ten\\
 &\dots\ten S(h_1)\sw{n-2}h_n\\
 &=S(h_1)\sw{1}h_2\ten\dots\ten \beta(\epsilon(S(h_1)
 \sw{i-1}\beta(\epsilon(h_i))))(S(h_1)\sw{i-2})h_{i-1}\ten \\
 &\dots\ten S(h_1)\sw{n-2}h_n\\
 &=S(h_1)\sw{1}h_2\ten\dots\ten \beta(\epsilon(S(h_1)
 \sw{i-1}(h_i)))(S(h_1)\sw{i-2})h_{i-1}\ten \\
 &\dots\ten S(h_1)\sw{n-2}h_n=\sigma_{i-1}\tau_n(h_1\ten\dots\ten h_n).
 \end{align*}
 The converse   is evident because
 $$ T(h\sw{1})\sw{1}h\sw{2}\ten T(h\sw{1})\sw{2}=\tau_2^3(1_H\ten T(h)).$$

\end{proof}

Let $R$ be a $k$-algebra and $\mathcal{H}$ a $k$-Hopf algebra endowed
with a twisted antipode $\widetilde{S}_\delta$ such that
${{\widetilde{S}}}_\delta^2=i d_{\mathcal{H}}$.
We  have a para-Hopf algebroid structure on $R \otimes \mathcal{H}\otimes R^{op}$ as in
 Example \ref{reza}. In  the next proposition  we recall the computation of the  cyclic cohomology of $R\otimes
 \mathcal{H}\otimes R^{op}$ as a para-Hopf algebroid in terms of   Hopf-cyclic cohomology of  $\mathcal{H}$.

\begin{prop}(\cite{achm01})\label{reza2}
Let $\mathcal{H}$ be a Hopf algebra as above. Then $$ {HC^*}(R\otimes\mathcal{H}\otimes
R^{op})=HC_{(\delta,1)}^*(\mathcal{H}).$$
\end{prop}

\begin{definition}{(Haar system for bialgebroids).} Let $(H,R)$ be a bialgebroid.
Let $\tau : H\longrightarrow R$ be a  right R-module map. We call $\tau$ a left
 Haar system for $H$ if  for all $h\in H$
$$\alpha(\tau(h^{(1)})) h^{(2)}=\alpha(\tau(h))1_H,$$
and $\alpha \tau=\beta\tau$.
We call $\tau$ a normal left Haar system if $\tau(1_H)=1_R.$
\end{definition}

We give a few examples of Haar systems. Let $H$  be the extended
Hopf algebra of a groupoid with finite base (Example \ref{reza3}).
Then it is easy to see that $\tau :H\rightarrow R $
 defined by $\tau(id_x)=id_x$  for
all $x\in Obj(\mathcal{G})$ and $\tau (\gamma)=0$ if $\gamma$ is
not an identity morphism, is a normal Haar system for $H$. For a
second example,
  one can directly  check that  the map
  $\tau : A_\theta \rightarrow \mathbb{C}\lbrack
 U,U^{-1}\rbrack$
defined by $$\tau(U^nV^m)=\delta_{m,0}U^n$$ is a normal Haar system for
the noncommutative torus $A_\theta $.
 \begin{prop}
Let  $H$ be a para-Hopf algebroid that admits a normal left Haar system. Then
$HC^{2i+1}(H)=0$ and $HC^{2i}(H)=\ker(\alpha-\beta)$ for all $i\geq 0$.
\end{prop}
\begin{proof}
Let $\eta :H^{\otimes_Rn}\longrightarrow H^{\otimes_R(n-1)}$ be the map
  $$\eta (h_1\otimes_R \dots \otimes_R h_n)=
\alpha(\tau(h_1)h_2\otimes_R \dots \otimes_R h_n.$$
 It is easy to check that $\eta $ is a contracting homotopy for
the Hochschild complex of $H_{\natural}$ and hence $HH^n({H})=0$
for $n>0$ and $HH^0({H})=\ker(\alpha-\beta)$. The rest  follows
from Connes's long exact sequence relating Hochschild and cyclic
cohomology.
\end{proof}
\begin{corollary}
Let $H$ be the para-Hopf algebroid of Example \ref{reza3}.
Then  $HC^{2i+1}(H)=0$ and $HC^{2i}(H)=R$ for
all $i\geq 0$.
\begin{multline*}
~\square
\end{multline*}
\end{corollary}
\begin{corollary}
We have
\begin{center}
$HC^{2i+1}(A_{\theta})=0$ and $HC^{2i}(A_{\theta})=\mathbb{C}\lbrack U,U^{-1}\rbrack$ for all $i\geq 0$.
\end{center}
\vspace{-1cm}
\begin{multline*}
~\square
\end{multline*}
 \end{corollary}

\end{document}